\documentclass[12pt,reqno]{amsart}
\usepackage{amssymb,amsthm,amsfonts,amstext}
\usepackage{amsmath}
\usepackage{zito}
\usepackage{mathrsfs}
\usepackage{fourier}

\usepackage{dsfont}

\usepackage{color}

\makeatletter \@addtoreset{equation}{section}
\@addtoreset{figure}{section}
 \makeatother

\newtheorem{theorem}{Theorem}[section]
\newtheorem{lemma}[theorem]{Lemma}
\newtheorem{proposition}[theorem]{Proposition}

\theoremstyle{remark}
\newtheorem{remark}[theorem]{Remark}

\newcommand{\mc}[1]{{\mathcal #1}}

\newcommand{\bb}[1]{{\mathbb #1}}

\newcommand{\<}{\langle}
\renewcommand{\>}{\rangle}

\newcommand{\cadlag}{{c\`adl\`ag~}}

\newcommand{\tdn}{{\mathbb{T}_n^d}}
\newcommand{\onk}{{\Omega_{n,k}}}
\newcommand{\sumxy}{{\sum_{\substack{x, y \in \tdn \\ y \neq x}}}}
\newcommand{\partialt}{{\partial_{\!t}}}
\newcommand{\norm}[1]{{\| #1 \|}}
\newcommand{\tnorm}[1]{{|\hspace{-1.5pt}|\hspace{-1.5pt}| #1 |\hspace{-1.5pt}|\hspace{-1.5pt}|}}

\title{Quantitative hydrodynamics for a generalized contact model}

\author{Julian Amorim}
\address{Instituto Nacional de Matem\'atica Pura e Aplicada, IMPA, Estrada Dona Castorina 110, 22460-320  Rio de Janeiro, Brazil.}  \email{julian.alexandre@impa.br}

\author{Milton Jara}
\address{Instituto Nacional de Matem\'atica Pura e Aplicada, IMPA, Estrada Dona Castorina 110, 22460-320  Rio de Janeiro, Brazil.}  \email{mjara@impa.br}

\author{Yangrui Xiang}
\address{Instituto Nacional de Matem\'atica Pura e Aplicada, IMPA, Estrada Dona Castorina 110, 22460-320  Rio de Janeiro, Brazil.}  \email{yangrui.xiang@impa.br}

\begin{document}

\begin{abstract}
We derive a quantitative version of the hydrodynamic limit obtained in \cite{ChaD-MLebPre} for an interacting particle system  inspired by integrate-and-fire neuron models.
More precisely, we show that the $L^2$-speed of convergence of the empirical density of states in a generalized contact process defined over a $d$-dimensional torus of size $n$ is of the optimal order $\mc O(n^{d/2})$. In addition, we show that the typical fluctuations around the aforementioned hydrodynamic limit are Gaussian, and governed by a inhomogeneous stochastic linear equation.
\end{abstract}

\maketitle

\section{Introduction}

A fundamental problem in non-equilibrium statistical mechanics is the derivation of macroscopic evolution equations from underlying microscopic dynamics. The aim is to identify effective variables that describe the large-scale behavior of many-particle systems and to characterize the evolution of such variables. Since the goal is the characterization of the evolution of the effective variables, a theoretical simplification, which makes the problem mathematically tractable, is to assume that the microscopic evolution is stochastic. in this case the mathematical framework of \emph{hydrodynamic limits of interacting stochastic systems}, see \cite{KipLan} for a comprehensive review, allows to obtain rigorous results into this regard. In a recent development \cite{JarMen}, a theory of \emph{quantitative hydrodynamics} has been proposed, and our aim in this article is to apply quantitative hydrodynamics to a stochastic system inspired by neuronal integrate-and-fire models, see Chapter 4 of \cite{GerKis}. This model, called by the authors a \emph{generalized contact process} has been introduced in \cite{ChaLeb} and the corresponding hydrodynamic limit was derived in \cite{ChaD-MLebPre}. Let us briefly describe the dynamics of this model; rigorous definitions are found in Section \ref{s2.1}. Neurons are represented by elements of a subset $\Lambda$ of the discrete lattice $\bb Z^d$. To each neuron $x$ we associate an electric potential $\sigma_x$, which for simplicity assumes integer values between $0$ and a threshold $k$. Neurons with potential $k$ are considered to be \emph{active}. Active neurons interact with other neurons, emitting a signal that raises the potential of the receiving neuron. The interaction between neurons at sites $x,y$ is modulated by a kernel $J_{x,y}$. Active neurons become inactive by losing all of its electric potential after an exponential time of rate $a$. 

In \cite{ChaLeb}, \cite{ChaD-MLebPre}, the authors work in infinite volume. In order to obtain quantitative results, it is more convenient to work in finite volume, so our setting will differ a bit from the setting of \cite{ChaLeb}, \cite{ChaD-MLebPre}. 
We introduce a scaling parameter $n$ that we will send to $\infty$, that measures the size of the system. For simplicity, we will take $\Lambda$ as a box of size $n$ and we consider periodic boundary conditions. Notice, however, that the \emph{geometric information} of the model is contained in the interaction kernel $J$, which does not need to be translation invariant.

Our main result states a quantitative version of the hydrodynamic limit derived in \cite{ChaD-MLebPre}. We will show that the empirical density of states (which is a stochastic process) is well approximated by the solution of the non-linear evolution equation \eqref{echid} (which is a deterministic evolution equation). Actually, we can say more. If we look at any box of size increasing with $n$ (like, say $n^\delta$ with $\delta >0$), then the empirical density of states is close to the value of the solution of \eqref{echid} at the corresponding space-time point.

Another important problem in the theory of hydrodynamic limits, is the derivation of the scaling limit of the \emph{fluctuations around its hydrodynamic limit} of the empirical density of states. When speaking about fluctuations, one can think about \emph{large fluctuations}, which corresponds to a large deviations principle, or one can think about \emph{typical fluctuations}, in which case we are dealing with a central limit theorem. In this article, we study the latter, and we prove that the typical fluctuations of the model follow a linear SPDE, which corresponds to the linearization of the hydrodynamic equation, perturbed by a space-time white noise of variance given in terms of the solution of \eqref{echid}. A general theory of typical fluctuations was a longstanding open problem in the theory of hydrodynamic limits, and this turns out to be the main application for which quantitative hydrodynamics was introduced. 

This article is organized as follows. In Section \ref{s2} we introduce the model considered here in a rigorous way, as well as we state our main results in a mathematically precise way. In Theorem \ref{t1}, we state our quantitative version of the hydrodynamic limit for the empirical density of states. We point out that the quantitative estimate of Theorem \ref{t1} is not sharp in $n$, due to the fact that the size of the error of a relevant discretization of the hydrodynamic equation is larger than the size of typical fluctuations. In order to obtain an estimate with an optimal error, we need to consider the fluctuations around the discretization of the hydrodynamic equation, which \emph{a fortiori} seems to be the right thing to look at. This estimate is stated in Theorem \ref{t1b}. 

We see in Theorem \ref{t1} that the time dependence of our quantitative estimate is given by an exponential tower $e^{Ce^{Ct}}$. This is somehow the weakest point of our methodology. The idea is that the model has a plethora of absorbing states; any initial state without active sites is an absorbing state. This is reflected in the hydrodynamic equation, which also has states without active sites as absorbing states. In \cite{ChaLeb} the authors analyze the question of long-term behavior of solutions of \eqref{echid}, and they show that for some choices of the parameters of the model, \eqref{echid} has non-trivial steady states, that is, states with non-zero active states. Since our methodology in principle works for arbitrary interactions between states, one can not distinguish between hydrodynamic equations with or without non-trivial steady states. Observe that without further assumptions on the hydrodynamic equation, the best one can do is to show that solutions approach the trivial steady states at most at exponential speed. For small state densities, fluctuations transit from Gaussian to Poissonian, and in particular we no longer have concentration estimates to our disposal. An interesting line of research would be the analysis of long-term behavior of the model, which are out of range to the techniques used in this article.

Theorem \ref{t2} states the convergence of the typical fluctuations to solutions of the corresponding SPDE. Since this SPDE is linear and has Gaussian noise, the typical fluctuations are given by a space-time Gaussian process, whose correlation structure depends on the solution to the hydrodynamic equation \eqref{echid}.

In Section \ref{s3} we derive quantitative hydrodynamics for the model, that is, we prove Theorems \ref{t1} and \ref{t1b}. The proof follows more or less immediately from Theorem \ref{t3}, which is the main innovation of \cite{JarMen}. Theorem \ref{t3} compares the law of the model at times $t \geq 0$ with a non-homogeneous Bernoulli product measure, which is much simpler to analyze. This comparison is achieved using Yau's \emph{relative entropy method} introduced in \cite{Yau}, see Chapter 6 of \cite{KipLan}. Here we use the approach of \cite{JarMen}, which yield better bounds. 

In Section \ref{s4}, we prove Theorem \ref{t2}, that is, the central limit theorem for the fluctuations of the empirical measure around its hydrodynamic limit. The idea of the proof is to show that the fluctuations of the empirical density of particles satisfy a discrete approximation of the limiting SPDE \eqref{SPDE}. The various error terms are shown to be negligible using the entropy estimate of Theorem \ref{t3} as an \emph{a priori} estimate.

In Appendix \ref{a1} we review the solution theory of the hydrodynamic equation \eqref{echid}. One can understand \eqref{echid} as a non-linear, non-local  transport equation, and therefore some care is needed in order to state existence, uniqueness and regularity properties of solutions of \eqref{echid}, as well as convergence properties of suitable discretizations. In consequence, we have decided to include complete proofs instead of referring to classical works which may not be immediately applicable to our situation without adaptation of some arguments. In Appendix \ref{apb} we collect some well known results about relative entropy and concentration of averages of independent random variables, which are used throughout the text.

\section{Preliminaries and Notation}
\label{s2}
\subsection{The Generalized Contact Process}
\label{s2.1}
Let $n \in \bb N$ be a scaling parameter. Let $\tdn := \bb Z^d/ n \bb Z^d$ be the discrete torus of size $n$. The elements $x \in \tdn$ are called \emph{sites} and represent the location of either neuron on a integrate-and-fire model or individuals in a population model. Let $k \in \bb N$ be a fixed parameter and let $Q_k:= \bb Z / (k+1) \bb Z$ the space of integers modulo $k+1$. We will identify $Q_k$ with the set $\{0,1,\dots,k\}$. Let $\onk := Q_k^{\tdn}$ and denote by $\sigma := (\sigma_x; x \in \tdn)$ the elements of $\onk$. Observe that $\onk$ has a ring structure. The value of $\sigma_x$ represents the voltage of the neuron at site $x$, or the level of infection of an individual at site $x$. We say that site $x$ is \emph{active} if $\sigma_x = k$, otherwise we say that the site $x$ is \emph{passive}. Let $(\delta_x; x \in \tdn)$ be the canonical basis of the free module $\onk$, that is,
\[
(\delta_x)_y :=
\left\{
\begin{array}{r@{\;;\;}l}
1 & y =x \\
0 & y \neq x.
\end{array}
\right.
\]
Let $\bb T^d : = \bb R^d / \bb Z^d$ be the continuous torus. We consider $\frac{1}{n} \tdn$ as embedded into $\bb T^d$. Let $J: \bb T^d \times \bb T^d \to [0,\infty)$ be bounded. Let $a>0$ be fixed and for each $x \in \tdn$, let $c_x = c_x^{n,k}:\onk \to [0,\infty)$ be given by
\[
c_x(\sigma) := a \mathds{1}(\sigma_x=k) + \tfrac{1}{n^d} \sum_{y \neq x} J \big( \tfrac{\vphantom{y}x}{n} , \tfrac{y}{n}\big) \mathds{1}(\sigma_y =k,\sigma_x \neq k)
\]
for every $\sigma \in \onk$. For each $f: \onk \to \bb R$, let $L_n f: \onk \to \bb R$ be defined as
\[
L_n f(\sigma) := \sum_{x \in \tdn} c_x(\sigma) \big( f(\sigma+\delta_x) -f(\sigma)\big)
\]
for every $\sigma \in \onk$. The operator $L_n$ defined in this way turns out to be the generator of a Markov chain $(\sigma(t); t \geq 0)$ introduced in \cite{ChaLeb}, which the authors call the \emph{generalized contact process}. 

Observe that the set $\{ \sigma \in \onk; \sigma_x \neq k \text{ for every } x \in \tdn\}$ is the maximal absorbing set of the dynamics generated by $L_n$. Observe as well that for every initial datum $\sigma \in \onk$, the dynamics eventually arrives to this absorbing set. Our aim is to analyze the evolution of the chain in a time scale which in one hand is large enough to give rise to non-trivial evolutions of the empirical density of states, and in the other hand is short enough to avoid instantaneous absorption in the limit. Due to the non-linear nature and the lack of a smoothing mechanism of the hydrodynamic equation, the estimates of Theorem \ref{t1}, \ref{t1b} and \ref{t3}, although sharp in the scaling parameter $n$, deteriorate quite fast in time. In \cite{ChaLeb}, the authors analyze the presence of non-trivial stationary states for the hydrodynamic equation \eqref{echid}. An interesting problem which is out of reach for the moment, is the presence in the stochastic system of \emph{metastable} states corresponding to non-trivial states of the limiting equation.

Given a finite time horizon $T >0$, we denote by $\mc D([0,T]; \Omega_{n,k})$ the space of \cadlag trajectories with values in $\Omega_{n,k}$, endowed with the $J_1$-Skorohod topology and its associated Borel $\sigma$-algebra. For a given initial measure $\mu$ in $\Omega_{n,k}$, we denote by $\bb P_{\mu}^n$ the law in $\mc D([0,T]; \Omega_{n,k})$ of the process $(\sigma(t); t \geq 0)$ with initial measure $\mu$, and we denote by $\bb E_{\mu}^n$ the expectation with respect to $\bb P_{\mu}^n$.

\subsection{The hydrodynamic equation}

For $f: \bb T^d \to \bb R$ integrable, let $J \ast f: \bb T^d \to \bb R$ be given by
\begin{equation}
J \ast f(x) := \int J(x,y) f(y) dy
\end{equation}
for every $x \in \bb T^d$. Let $M$ be the $(k+1) \times (k+1)$-matrix given by
\[
M_{ij}:=
\left\{
\begin{array}{r @{\;;\;}l}
1 & i=1,\dots,k, j=i-1\\
-1 & i=0,\dots,k-1, j=i\\
0 & \text{otherwise.}
\end{array}
\right.
\]
Let $A$ be the $(k+1)\times (k+1)$-matrix given by
\[
A_{ij} := 
\left\{
\begin{array}{r@{\;;\;}l}
a & i=0, j=k\\
-a & i=k, j=k\\
0 & \text{otherwise.}
\end{array}
\right.
\]
The \emph{hydrodynamic equation} associated to the generalized contact model is the equation
\begin{equation}
\label{echid}
\partial_t \vec{u}(t,x) = A  \vec u(t,x) + \big(J \ast u^k(t,x)\big) M \vec u(t,x).
\end{equation}

In  \cite{ChaD-MLebPre} is was shown that the empirical density of sites with given state in the generalized contact process converge to a solution of \eqref{echid}. Our aim is to show a quantitative version of this result, as well as the corresponding central limit theorem. Let us state our main result in a rigorous way. Let $\mc P_k$ denote the space of probability measures in $Q_k$. Observe that $\mc P_k$ can be identified with a simplex in $\bb R^{Q_k}$. For a given function $\vec u: \tdn \to \mc P_k$, we define the \emph{profile measure} $\mu_{\vec u}^n$ associated to $\vec u$ as
\[
\mu_{\vec u}^n(\sigma) := \prod_{x \in \tdn} u_x^{\sigma_x}
\]
for every $\sigma \in \onk$.
Let $\vec u_0: \bb T^d \to \mc P_k$ be bounded. For each $n \in \bb N$, let $\mu_{\vec u_0}^n$ be the probability measure in $\onk$ given by
\[
\mu_{\vec u_0}^n (\sigma) := \prod_{x \in \tdn} u^{\sigma_{\!\!x}}_0\big(\tfrac{x}{n} \big)
\]
for every $\sigma \in \onk$. For each $n \in \bb N$, we will consider the process $(\eta(t);t \geq 0)$ in $\onk$ with initial condition $\mu_{\vec u_0}^n$. In that case we say that $\vec u_0$ is the \emph{initial profile} of the processes $(\eta(t); t \geq 0)$. Observe that $\mu_{\vec u_0}^n$ coincides with the profile measure associated to the function $x \mapsto u_0(\frac{x}{n})$.

We will make the following assumptions on the initial profile $\vec u_0$:

\begin{equation}
\tag{H1}
\label{H1}
\vec{u}_0 \text{ is of class } \mc C^2 \text{ and there exists } \epsilon >0 \text{ such that } \vec u_0(x) \in [\epsilon,1-\epsilon]^{Q_k} \text{ for every } x \in \bb T^d.
\end{equation}
\\
We will prove the following result:

\begin{theorem}
\label{t1} There exists a constant $C = C(\epsilon)$ such that for every continuous function $f: \bb T^d \to \bb R$, every $i \in Q_k$ and every $t \geq 0$,
\[
\bb E_{\mu_{\vec{u}_0}^n}^n\Big[ \Big( \tfrac{1}{n^d} \!\!\!\sum_{x \in \tdn}\!\! \big( \mathds{1}\big(\eta_x(t)=i\big) - u^i\big(t, \tfrac{x}{n} \big) \big)f \big( \tfrac{x}{n} \big) \Big)^ 2 \Big] \leq \frac{Ce^{Ce^{Ct}} \|f\|_\infty^2}{n^{\min\{2,d\}}},
\]
where $(\vec u(t,x); t \geq 0, x \in \bb T^d)$ is the solution of \eqref{echid} with initial condition $\vec u_0$.
\end{theorem}

This result can be understood as a law of large numbers for the empirical density. The proof of this result is based on the so-called \emph{Yau's relative entropy method}, which aims to estimate the relative entropy between the law of the process at time $t$ and properly defined profile measures. 

We point out that the dependence in $n$ of the estimate in Theorem \ref{t1} is not optimal; the optimal speed of convergence would be $n^d$. We will see in Theorem \ref{t1b} below that this is just a finite-size effect that can be overcome introducing better discrete approximations for the expectation of $\sigma_x(t)$.

The exponential tower in the estimate of Theorem \ref{t1} comes from the fact that \eqref{echid} has a quadratic nonlinearity and that the kernel $J$ does not have an evident smoothing effect in the equation. 

Once one has obtained a law of large number for a family of random variables, it is natural to ask about the fluctuations around this law of large numbers. If one is interested in typical fluctuations, one wants to discuss about the \emph{central limit theorem}. If one is, otherwise, interested in atypical fluctuations, then the objective is to derive a \emph{large deviations principle}. In this article, we will aim to derive a central limit theorem for the empirical voltage of neurons.

Fluctuating hydrodynamics proposes that the typical fluctuations around its hydrodynamic limit of an interacting particle system are governed by a stochastic PDE obtained by linearizing the hydrodynamic equation and adding a Gaussian noise:
\begin{equation}
\label{SPDE}
d \!\vec X_t = \Big( A \vec X_t + \big(J \ast u^k_t\big) M \vec X_t + \big(J \ast X^k_t \big) M \vec u_t \Big) d\!t + d\hspace{-1pt} \vec{\mc M}_t,
\end{equation}
where $(\vec{\mc M}_t; t \geq 0)$ is a Gaussian martingale. 
For each $t \geq 0$ and $x \in \bb T^d$, let $\gamma(t,x)$ be the $(k+1) \times (k+1)$ symmetric matrix such that
\[
\gamma^{i,j}(t,x) := 0 \text{ whenever } j-i \neq -1,0,1,
\]
that is, $\gamma$ is \emph{tridiagonal},
\[
\gamma^{i,i}(t,x) :=
\left\{
\begin{array}{r@{\;;\;}l}
a u^k(t,x) +  \big (J \ast u^k(t,x) \big) u^{k-1}(t,x) & i=k\\
a u^k(t,x) +  \big (J \ast u^k(t,x) \big) u^{0}(t,x) & i=0\\
\big (J \ast u^k(t,x) \big) \big( u^{i-1}(t,x) + u^i(t,x)\big) & i \neq 0,k,
\end{array}
\right.
\]
and
\[
\gamma^{i,i+1}(t,x) = \gamma^{i+1,i}(t,x):= 
\left\{
\begin{array}{c@{\;;\;}l}
- a u^k(t,x)   & i=k\\
-\big (J \ast u^k(t,x) \big) u^i(t,x) & i \neq k.
\end{array}
\right.
\]
The integrated quadratic covariation of the martingale $( \vec{\mc M}_t; t \geq 0)$ is given by
\begin{equation}
\label{Gauss}
\big\<\mc M^i\!(f), \mc M^j\!(g) \big\>_t := \int_0^t \int f(x) g(x) \gamma^{i,j}(s,x) dx ds
\end{equation}
for every $f,g \in \mc C(\bb T^d; \bb R)$ and every $t \geq 0$.

It will be convenient to center the empirical density not around $\vec u(t, \frac{x}{n})$, but around a properly defined approximation of the solution of \eqref{echid}.
Let us define $J^n : \tdn \times \tdn \to [0,\infty)$ as
\[
J^n_{x,y} := J\big( \tfrac{\vphantom{y}x}{n}, \tfrac{y}{n}\big) \text{ for every } x,y \in \tdn \text{ such that } x \neq y.
\]
For $f: \tdn \to \bb R$, let us define $J^n \ast f: \tdn \to \bb R$ as
\[
(J^n \ast f)_x := \frac{1}{n^d} \sum_{y \in \tdn} J_{x,y}^n f_y \text{ for every } x \in \tdn.
\]
Let $(\vec{u}_x^{n}(t); x \in \tdn, i \in Q_k, t \geq 0)$ be the solution of the equation
\begin{equation}
\label{numaprox}
\tfrac{d}{dt} \vec{u}_x^n(t) = A \vec{u}_x^n(t) + \big( J^n \ast u^{n,k}(t)\big)_x M \vec{u}_x^n(t)
\end{equation}
with initial condition $x \mapsto \vec{u}_0(\frac{x}{n})$. This equation serves as a numerical approximation of \eqref{echid}. 

If we center the empirical density of states around the functions $\vec u_t^n$, we obtain a version of Theorem \ref{t1} with better error:

\begin{theorem}
\label{t1b}
There exists a constant $C = C(\epsilon)$ such that for every continuous function $f: \bb T^d \to \bb R$, every $i \in Q_k$ and every $t \geq 0$,
\[
\bb E_{\mu_{\vec{u}_0}^n}^n\Big[ \Big( \tfrac{1}{n^d} \!\!\!\sum_{x \in \tdn}\!\! \big( \mathds{1}\big(\eta_x(t)=i\big) - u^{n,i}_x(t) \big) f \big( \tfrac{x}{n} \big) \Big)^ 2 \Big] \leq \frac{Ce^{Ce^{Ct}} \|f\|_\infty^2}{n^d},
\]
\end{theorem}

For each $t \geq 0$, each $i \in Q_k$ and $f \in \mc C^\infty(\bb T^d)$, let $X_t^{n,i}(f)$ be given by
\[
X_t^{n,i}\!(f) := \frac{1}{n^{d/2}} \!\!\!\sum_{x \in \tdn}\!\! \Big( \mathds{1}\big(\eta_x(t)=i\big) - u^{n,i}_x(t) \Big)f \big( \tfrac{x}{n} \big).
\]
This relation defines a distribution-valued process $(\vec{X}_t^n; t \geq 0)$ which we call the \emph{fluctuation field} associated to the empirical density. For $\vec f \in \mc C^\infty(\bb R^d; \bb R^{Q_k})$, it will be useful to introduce the notation
\[
\vec X_t^{n}(\vec f) := \sum_{i \in Q_k} X_t^{n,i} (f^i).
\]
Observe that the bound of Theorem \ref{t1} can be restated as
\[
\bb E_{\mu_{\vec{u}_0}^n}^n\big[ \big\| \vec X_t^{n} (\vec f) \big\|^2\big] \leq Ce^{Ce^{Ct}} \|\vec f\|_\infty^2.
\]
In other words, $\vec X_t^{n}(\vec f)$ is uniformly bounded in $L^2(\bb P_{\!\!\mu_{\vec{u}_0}^n}\!)$. Observe that under $\mu_{\vec{u}_0}^n$, $\vec{X}_0^n$ converges to a Gaussian process $\vec{\xi}$ of covariance given by
\begin{equation}
\label{cov1}
\bb E \big[ \xi^i(f) \xi^j(g) \big] = -\int f(x) g(x) u_0^i(x)u_0^j(x) dx
\end{equation}
for $i \neq j$ and
\begin{equation}
\label{cov2}
\bb E \big[ \xi^i(f) \xi^i(g) \big] = \int f(x) g(x) u_0^i(x)\big(1-u_0^i(x)\big) dx.
\end{equation}
We have the following result:

\begin{theorem}
\label{t2}
The sequence of processes $( \vec{X}_t^n; t \geq 0)_{n \in \bb N}$ converge in distribution of the solution of \eqref{SPDE} with initial condition $\vec \xi$, with respect to the $J_1$-Skorohod topology of paths in the space of distributions.
\end{theorem}

\section{Hydrodynamic Limit}
\label{s3}
As mentioned above, the proof of Theorem \ref{t1} is based on a quantitative estimate on the relative entropy of the law of $\sigma(t)$ and a properly defined profile measure. 

For each $t \geq 0$, let $\mu_t^n$ be the profile measure associated to $\vec{u}^n(t)$, that is,
\[
\mu_t^n(\sigma) := \prod_{x \in \tdn} u_x^{n,\sigma_x}(t) \text{ for every } \sigma \in \onk.
\]
Observe that $\mu_0^n = \mu_{\vec{u}_0}^n$. For each $t \geq 0$, let $f_t^n: \onk \to [0,\infty)$ be the density of the law of $\sigma(t)$ with respect to $\mu_t^n$, that is,
\[
f_t^n(\sigma) := \frac{\bb P_{\mu_0^n}^n \big( \sigma(t)=\sigma \big)}{\mu_t^n(\sigma)} \text{ for every } \sigma \in \onk.
\]
Let $H_n(t)$ be the relative entropy of $f_t^n$, that is,
\[
H_n(t) := \int f_t^n \log f_t^n d \mu_t^n \text{ for every } t \geq 0.
\]
We will prove the following estimate:

\begin{theorem}
\label{t3} 
There exists a finite constant $C = C(\epsilon, a, J)$ such that
\[
H_n(t) \leq C \big(e^{C(e^{Ct}-1)}-1 \big) \text{ for every } t \geq 0.
\]
\end{theorem}

\subsection{Yau's relative entropy method} In this section, we will prove Theorem \ref{t3}. The idea is to use \emph{Yau's relative entropy method}, introduced by H.T.~Yau in \cite{Yau}, with some elements taken from \cite{JarMen}. We briefly recall Yau's entropy inequality:

\begin{proposition}[Yau's relative entropy inequality]
\label{p1}
We have that
\begin{equation}
\label{coquimbo}
\tfrac{d}{dt} H_n(t) \leq \int \big( L^\ast \mathds{1} - \partial_{\!t} \log \mu_t^n \big) f_t^n d \mu_t^n,
\end{equation}
where $L^\ast$ is the adjoint operator of $L$ in $L^2(\mu_t^n)$.
\end{proposition}

This lemma corresponds to Lemma 6.1.4 in \cite{KipLan}. Our aim is to compute
\[
\bb F_{\!t}:= L^\ast \mathds{1} - \partial_{\!t} \log \mu_t^n.
\]
and to obtain a bound on $\int \bb F_{\!t} d \mu_t^n$ in terms of $H_n(t)$. First we need to compute $\bb F_{\!t}$ explicitly. In order to simplify the notations, let us use the notation $u_x^i := u_x^{n,i}(t)$ and let us define
\begin{equation}
\label{laserena}
w_x^i:= \mathds{1}(\sigma_x=i) -u_x^{i}
\end{equation}
for every $x \in \tdn$, every $i \in Q_k$ and every $t \geq 0$. 

\begin{lemma}
\label{l1}
For every $t \geq 0$, we have that
\[
\bb F_t =  \frac{1}{n^d}  \sum_{i=0}^{k} \sumxy J^n_{x,y} g_x^i w_x^i w_y^k, 
\]
where
\[
g_x^i:=
\left\{
\begin{array}{c@{\;;\;}l}
\frac{u_x^{k-1}}{u_x^k} & i=k\\
-1 & i=0\\
\frac{u_x^{i-1}-u_x^i}{u_x^i} & i \neq 0,k.
\end{array}
\right.
\]
\end{lemma}

\begin{proof}
First we observe that the inverse of the transformation $\sigma \mapsto \sigma+\delta_x$ is the transformation $\sigma \mapsto \sigma-\delta_x$. Therefore, for every measure $\mu$ in $\onk$, 
\[
L^\ast \mathds{1} = \sum_{x \in \tdn} \Big(c_x(\sigma-\delta_x) \frac{\mu(\sigma-\delta_x)}{\mu(\sigma)} -c_x(\sigma)\Big).
\]
Now we observe that
\[
\frac{\mu_t^n(\sigma-\delta_x)}{\mu_t^n(\sigma)} = \frac{u_x^{\sigma_x-1}}{u_x^{\sigma_x}}.
\]
Therefore, 
\[
\begin{split}
L^\ast \mathds{1} 
	&= \sum_{x \in \tdn}  a \Big( \frac{\mathds{1}(\sigma_x=0)}{u_x^0} -\frac{\mathds{1}(\sigma_x=k)}{u_x^k}\Big) u_x^k\\
	&\quad \quad + \frac{1}{n^d} \sumxy \Big( \frac{\mathds{1}(\sigma_x=i)}{u_x^i} - \frac{\mathds{1}(\sigma_x=i-1)}{u_x^{i-1}}\Big) u_x^{i-1} J_{x,y}^n \mathds{1}(\sigma_y = k).
\end{split}
\]
Similarly,
\[
\begin{split}
\partial_{\!t} \log \mu_t^n 
	&= \sum_{x \in \tdn} \frac{\partial_{\!t} u_x^{\sigma_x}}{u_x^{\sigma_x}}\\
	&= \sum_{x \in \tdn} a u^k_x \Big(  \frac{\mathds{1}(\sigma_x=0)}{u_x^0} -\frac{\mathds{1}(\sigma_x=k)}{u_x^k}\Big) \\
	&\quad \quad + \frac{1}{n^d} \sum_{i=1}^k \sumxy \Big( \frac{\mathds{1}(\sigma_x=i)}{u_x^i} - \frac{\mathds{1}(\sigma_x=i-1)}{u_x^{i-1}}\Big) u_x^{i-1} J_{x,y}^n u_y^k.
\end{split}
\]
Observe now that
\[
\frac{\mathds{1}(\sigma_x=i)}{u_x^i} - \frac{\mathds{1}(\sigma_x=i-1)}{u_x^{i-1}} = \frac{w_x^i}{u_x^i} -\frac{w_y^{i-1}}{u_x^{i-1}}.
\]
We conclude that
\[
\bb F_{\!t} 
		= \frac{1}{n^d} \sum_{i=1}^k \sumxy 
				u_x^{i-1} J_{x,y}^n \Big( \frac{w_x^i}{u_x^i} -\frac{w_y^{i-1}}{u_x^{i-1}}\Big) w_y^k,
\]
as we wanted to show.
\end{proof}

The main consequence of this result is that we can write $\bb F_{\!t}$ as a sum of monomials of degree two with respect to the variables $(w_x^i; x \in \tdn, i \in Q_k)$. We will see below that this decomposition in monomials of degree two plays a fundamental role in the proof of Theorem \ref{t3}.

\subsection{Proof of Theorem \ref{t3}}

With the aid of Proposition \ref{p1} and Lemma \ref{l1}, we can prove Theorem \ref{t3}. By Lemma \ref{l1} and Yau's estimate \eqref{coquimbo},
\[
\tfrac{d}{dt} H_n(t) \leq \int \bb F_t f_t d \mu_t.
\]
It will be convenient to introduce the notation
\[
\bb F_t^i :=  \frac{1}{n^d} \sumxy J^n_{x,y} g_x^i w_x^i w_y^k, 
\]
By the entropy inequality, for every $\gamma >0$,
\[
\int \bb F_t^i f_t d \mu_t \leq \gamma^{-1} \Big( H_n(t) + \log \int e^{\gamma \bb F_t^i} d \mu_t \Big).
\]
Under $\mu_t$, the random variables $((w_x^i, w_x^k); x \in \tdn)$ are independent. Therefore, the integral $\int e^{\gamma \bb F_t^i} d \mu_t $ can be estimated by Hanson-Wright inequality (see Proposition \ref{pa2}). First we observe that $w_x^i$ is a centered Rademacher variable, and therefore by Hoeffding's inequality it is sub-Gaussian of index $\frac{1}{4}$. Therefore, for 
\[
\gamma^{-1} \geq \Big(\frac{1024}{16n^{2d}} \sumxy \big( J_{x,y}^n\big)^2 (g_x^i)^2\Big)^{1/2},
\]
\[
\log \int e^{\gamma \bb F_t^i} d \mu_t \leq \log 3.
\]
By Lemma \ref{la4} (see also the discussion after estimates \eqref{talca}, \eqref{linares}), $g_x^i \leq \epsilon_0^{-1} e^{K_nt}$, where $K_n:= \max\{a,\norm{J}_{1,n}\}$. Therefore,
\[
\begin{split}
\frac{1024}{16n^{2d}} \sumxy \big( J_{x,y}^n\big)^2 (g_x^i)^2 \leq 64 \epsilon_0^{-2} \norm{J}_\infty^2 e^{2K_nt},
\end{split}
\]
and therefore we can take $\gamma^{-1} = 8 \epsilon_0^{-1} \norm{J}_\infty e^{K_nt}$ to conclude that
\[
\tfrac{d}{dt} H_n(t) \leq 8(k+1)\epsilon_0^{-1} \norm{J}_\infty e^{K_nt}\big( H_n(t) + \log 3\big).
\]
Integrating this inequality, we conclude that
\[
H_n(t) \leq \big( H_n(0)+\log 3 \big) \exp \big\{ \tfrac{8(k+1)\epsilon_0^{-1} \norm{J}_\infty}{K_n}(e^{K_nt}-1) \big\} -\log 3.
\]
Since $K_n$ is uniformly bounded in $n$, we have proved Theorem \ref{t3}.

\subsection{Proof of  Theorem \ref{t1}}

In this section we will show that Theorem \ref{t1} is a simple consequence of Theorem \ref{t2}. By Lemma \ref{pitagoras} and then Lemma \ref{hoeff}, 
\[
\Big\| \frac{1}{n^d} \sum_{x \in \tdn} w_x^i f \big( \tfrac{x}{n}\big) \Big\|_{\psi_2}^2 \leq \frac{1}{n^{2d}} \sum_{x \in \tdn} \|w_x^i \|_{\psi_2}^2 f \big( \tfrac{x}{n}\big)^2 \leq \frac{1}{4n^{2d}} \sum_{x \in \tdn}  f \big( \tfrac{x}{n}\big)^2
\]
Using \eqref{donvar}, we see that
\[
\bb E_{\mu_{\vec{u}_0}^n}^n\Big[ \Big(\frac{1}{n^d} \sum_{x \in \tdn} w_x^i f \big( \tfrac{x}{n}\big) \Big)^2\Big]
		\leq \gamma^{-1}\Big( H_n(t) + \log \int \exp\Big\{ \gamma \Big(\frac{1}{n^d} \sum_{x \in \tdn} w_x^i f \big( \tfrac{x}{n}\big) \Big)^2 \Big\} d \mu_t^n \Big)
\]
By Lemma \ref{quad},
\[
 \log \int \exp\Big\{ \gamma\Big(\frac{1}{n^d} \sum_{x \in \tdn} w_x^i f \big( \tfrac{x}{n}\big) \Big)^2 \Big\} d \mu_t^n 
 		\leq \log 3
\]
for
\[
\gamma^{-1} = \frac{1}{n^{2d}} \sum_{x \in \tdn} f \big(\tfrac{x}{n}\big)^2,
\]
from where we conclude that
\[
\bb E_{\mu_{\vec{u}_0}^n}^n\Big[ \Big(\frac{1}{n^d} \sum_{x \in \tdn} w_x^i f \big( \tfrac{x}{n}\big) \Big)^2\Big]
		\leq \frac{1}{n^{2d}} \sum_{x \in \tdn} f \big(\tfrac{x}{n}\big)^2 \Big( C\big(e^{C(e^{Ct}-1)}-1\big)+\log 3\Big),
\]
which proves the lemma.

\begin{remark}
\label{r1}
Observe that we have proved Theorem \ref{t1} with the norm $\|f\|_\infty$ replaced by the stronger \emph{discrete $\ell^2$-norm}
\[
\|f\|^2_{\ell_n^2} := \frac{1}{n^d} \sum_{x \in \tdn} f \big( \tfrac{x}{n}\big)^2.
\]
We have decide to state Theorem \ref{t1} with the norm $\|f\|_\infty$ in order to express the dependence on $f$ in terms of a norm that does not depend on $n$.
\end{remark}

\section{Non-equilibrium fluctuations}
\label{s4}
In this section we prove Theorem \ref{t2}. 
Let $(\mc F_t;t \geq 0)$ be a filtration, let $(\vec {\mc M}_t; t \geq 0)$ a Gaussian martingale with respect to this filtration with covariations given by \eqref{Gauss} and let $\vec \xi$ be an $\mc F_0$-measurable Gaussian process with covariances given by \eqref{cov1} and \eqref{cov2}. We say that a distribution-valued process $(\vec X_t; t \geq 0)$ adapted to $(\mc F_t;t \geq 0)$ is a solution of \eqref{SPDE} with initial condition $\vec \xi$ if
\begin{itemize}
\item[a)] $\vec X_0$ has law $\vec \xi$,

\item[b)] for every $t \geq 0$ and every test function $\vec f \in \mc C^1([0,t]; \mc C(\bb T^d; \bb R^{Q_k}))$,
\[
\vec X_t( \vec f_t) = \vec X_0(\vec f_0) + \int_0^t \big\{\vec X_s\big( \partial_{\!s}\vec f_s + A^\ast \vec f_s + (J \ast u_s^k) M^\ast \vec f_s \big) + X_s^k( J^\ast \ast \< \vec u_s, M^\ast \vec f_s\>   \big\} ds 
		+ \int d \vec{ \mc M}_s(\vec f_s).
\]
\end{itemize}
Here $J^\ast: \bb T^d \times \bb T^d \to [0,\infty)$ is defined as $J^\ast(x,y) := J(y,x)$ for every $x,y \in \bb T^d$.

The idea is to show that $(\vec{X}_t^n; t \geq 0)$ is an approximated solution of \eqref{SPDE} with initial condition $\vec \xi$. In order to to that, for each test function $f \in \mc C^\infty(\bb T^d)$, we start by using Dynkin's formula to decompose the process $(\vec X_t^n(\vec f); t \geq 0)$ into a drift part and a martingale part.

\subsection{The martingale decomposition}

Let $\Gamma_n$ be the \emph{carr\'e du champ} operator associated to $L_n$, that is, for every $f,g: \onk \to \bb R$, let $\Gamma_n(f,g): \onk \to \bb R$ be given by
\[
\Gamma_n(f,g) := L_n(fg) - f L_n g - g L_n f.
\]

According to Dynkin's formula, for every $g: [0,T] \times \onk \to \bb R^{Q_k}$ smooth on the time variable, the process
\[
g(t,\sigma(t)) - g(0,\sigma(0)) - \int_0^t \big(\partial_{\!s} + L_n \big) g(s, \sigma(s)) ds
\]
is an $\bb R^d$-valued martingale of predictable quadratic covariation
\[
\Big(\int_0^t \Gamma_n \big(g^i(s,\sigma(s)), g^j(s,\sigma(s)) \big)ds \Big)_{i,j \in Q_k}.
\]
Taking $g = (X_t^{n,i}(f); i \in Q_k)$ for each $ f \in \mc C^\infty(\bb T^d)$, we define the  martingale $(\vec M_t^{n}; t \geq 0)$ as
\begin{equation}
\label{dynkin}
M_t^{n,i}(f) :=  X_t^{n,i}(f) - \vec X_0^{n,i}(f) - \int_0^t \big( \partial_{\!s} + L_n \big) X_s^{n,i}(f) ds
\end{equation}
for every $t \geq 0$. The quadratic covariation of this martingale is equal to
\[
\< M^{n,i}(f), M^{n,j}(g)\>_t = \int_0^t \Gamma_n \big( X_s^{n,i}(f), X_s^{n,j}(g) \big) ds
\]

Our aim is to write the processes $(\partialt + L_n) X_t^{n,i}(f^i)$ and $\Gamma_n \big( X_s^{n,i}(f^i), X_s^{n,j}(f^j) \big)$ as functions of the process $(\vec X_t^n; t \geq 0)$ plus error terms, on which case we could show that $(\vec X_t^n; t \geq 0)$  satisfies an approximate version of a closed equation, which will turn out to be \eqref{SPDE}. Let us start with the quadratic covariation.

\begin{lemma}
\label{l2}
Under \eqref{H1}, for every $f, g \in \mc C^2(\bb T^d)$,
\[
\bb E_{\mu_{\vec u_0}^n}^n \Big[ \Big| \Gamma_n \big( X_s^{n,i}(f), X_s^{n,j}(g) \big) - \int \gamma^{i,j}(t,x) f(x) g(x) dx \Big|\Big] \leq \frac{C e^{C e^{Ct}} \norm{f}_\infty^2}{n^{\min\{2,d\!/\!2\}}}.
\]
In particular,
\[
\lim_{n \to \infty} \< M^{n,i}(f), M^{n,j}(g)\>_t = \int_0^t \int \gamma^{i,j}(s,x) f(x)g(x) dx ds
\]
in $\mc L^1(\bb P_{\mu_{\vec u_0}^n}^n)$.
\end{lemma}

\begin{proof}
Observe that
\[
\Gamma_n (f,g) (\sigma) = \sum_{x \in \tdn} c_x(\sigma) \big( f(\sigma +\delta_x) - f(\sigma) \big) \big( g(\sigma + \delta_x ) - g(\sigma) \big).
\]
By the polarization identity, without loss of generality, we can take $f=g$. Therefore,
\begin{equation}
\label{vicuna}
\begin{split}
\Gamma_n \big( X_s^{n,i}(f), X_s^{n,j}(f) \big)
		&= \frac{1}{n^d} \!\!\sum_{x \in \tdn}\!\! c_x(\sigma) \big( \mathds{1}(\sigma_x \!=\! i\!-\!1) \!-\!\mathds{1}(\sigma_x =i) \big)
				\big( \mathds{1}(\sigma_x \!=\! j\!-\!1) \!-\!\mathds{1}(\sigma_x \!=\!j) \big) f \big( \tfrac{x}{n} \big)^2\\
		&\!\!\!\!\!\!\!\!= \frac{1}{n^d}\!\!\sum_{x \in \tdn}\!\! c_x(\sigma) \big( \mathds{1}(\sigma_x= i\!-\!1)(\delta_{i,j}- \delta_{i-1,j}) + \mathds{1}(\sigma_x =i) (\delta_{i,j}-\delta_{i,j-1}) \big) f \big( \tfrac{x}{n} \big)^2.\\
\end{split}
\end{equation}
Recall definition \eqref{laserena}. Define $(\gamma^{n,i,j}_x(t); x \in \tdn, i,j \in Q_k, t \geq 0)$ as 
\[
\gamma_x^{n,i,j}(t) :=0 \text{ whenever } j-i \neq -1,0,1,
\]
\[
\gamma^{n,i,i}_x(t) :=
\left\{
\begin{array}{r@{\;;\;}l}
a u_x^{k} + \big( J^n \ast u^{k}\big)_x u_x^{k-1} & i=k\\
a u_x^{k} + \big( J^n \ast u^{k}\big)_x u_x^{0} & i=0\\
\big( J^n \ast u^{k}\big)_x \big( u_x^{i-1} + u_x^{i} \big) & i \neq 0,k,
\end{array}
\right.
\]
and
\[
\gamma_x^{n,i,i+1}(t) = \gamma_x^{n,i+1,i}(t) :=
\left\{
\begin{array}{c@{\;;\;}l}
- a u_x^{k} & i=k\\
- \big( J^n \ast u_x^{k} \big)_x u_x^{i} & i \neq k.
\end{array}
\right.
\]
Recall that we are using the notation $u_x^i$ instead of the more correct notation $u_x^{n,i}(t)$. We will adopt the same convention for $\gamma^n$ and we will write $\gamma_x^{ij}$ instead of $\gamma_x^{n,i,j}(t)$. Observe that the expectation of $\Gamma_n \big( X_s^{n,i}(f), X_s^{n,j}(f) \big)$ with respect to $\mu_t^n$ is equal to
\[
\frac{1}{n^d} \sum_{x \in \tdn} \gamma_x^{i,j} f \big( \tfrac{x}{n} \big)^2.
\]
The idea is to prove that \emph{with respect to $\bb P_{\mu_{\vec u_0}^n}^n$}, $\Gamma_n \big( X_s^{n,i}(f), X_s^{n,j}(f) \big)$ concentrates around its mean with respect to $\mu_t^n$. We will only explain in detail the case $i \neq k$, $j = i+1$. The other cases can be dealt with in the same way. For each $\sigma \in \onk$ and each $i \in Q_k$, let us define $\sigma^i \in \{0,1\}^{\tdn}$ as 
\begin{equation}
\label{losvilos}
\sigma_x^i := \mathds{1}(\sigma_x = i).
\end{equation}
We have that
\begin{equation}
\label{illapel}
\begin{split}
\Gamma_n \big( X_s^{n,i}(f), X_s^{n,i+1}(f) \big)
		- \frac{1}{n^d} \sum_{x \in \tdn} \gamma_x^{i,i+1} f \big( \tfrac{x}{n} \big)^2
			&= -\frac{1}{n^d}\sum_{x \in \tdn} \Big( \big( J^n \ast \sigma^k \big)_x \sigma_x^i - \big( J^n \ast u^k\big)_x u_x^i \Big) f \big( \tfrac{x}{n} \big)^2\\
			&= - \frac{1}{n^d} \sum_{x \in \tdn} \Big(\big( J^n \ast w^k \big)_x w_x^i + \big( J^n \ast u^k \big)_x w_x^i\\
			&\quad \quad \quad + \big( J^n \ast w^k \big)_x u_x^i \Big) f \big( \tfrac{x}{n} \big)^2.
\end{split}
\end{equation}
Each of the three terms in this sum can be estimated using the concentration inequalities of Appendix \ref{apb}. By \eqref{donvar2}, for every $\gamma >0$,
\[
\begin{split}
\bb E_{\mu_{\vec u_0}^n}^n \Big[ \Big| \frac{1}{n^d} \sum_{x \in \tdn} \big( J^n \ast w^k \big)_x w_x^i f \big( \tfrac{x}{n} \big)^2 \Big| \Big]
		&\leq \gamma^{-1}\Big( H_n(t) + \log 2 +\\
		&+ \max_{\pm} \;\log \int \exp \Big\{ \frac{\pm\gamma}{n^d} \sum_{x \in \tdn}  \big( J^n \ast w^k \big)_x w_x^i f \big( \tfrac{x}{n} \big)^2 \Big\} d \mu_t^n \Big).
\end{split}
\]
By Hanson-Wright inequality \eqref{pa2},
\[
\max_{\pm} \;\log \int \exp \Big\{ \frac{\pm\gamma}{n^d} \sum_{x \in \tdn}  \big( J^n \ast w^k \big)_x w_x^i f \big( \tfrac{x}{n} \big)^2 \Big\} d \mu_t^n \leq \log 3
\]
for
\[
\gamma^{-1} = 16 \Big(\frac{1}{n^{4d}} \sumxy \big( J_{x,y}^n)^2 f \big( \tfrac{x}{n} \big)^4 \Big)^{1/2} \leq \frac{16 \norm{J}_\infty \norm{f}_\infty^2}{n^d},
\]
from where we conclude, using Theorem \ref{t3}, that
\[
\bb E_{\mu_{\vec u_0}^n}^n \Big[ \Big| \frac{1}{n^d} \sum_{x \in \tdn} \big( J^n \ast w^k \big)_x w_x^i f \big( \tfrac{x}{n} \big)^2 \Big| \Big]
		\leq \frac{Ce^{Ce^{Ct}} \norm{f}_\infty^2}{n^d}.
\]
In a similar fashion,
\[
\begin{split}
\bb E_{\mu_{\vec u_0}^n}^n \Big[ \Big| \frac{1}{n^d} \sum_{x \in \tdn} \big( J^n \ast u^k \big)_x w_x^i f \big( \tfrac{x}{n} \big)^2 \Big| \Big] 
		&\leq \gamma^{-1} \Big( H_n(t) + \log 2 +\\
		&\quad \quad + \max_{\pm} \; \log \int \exp \Big\{ \frac{\pm\gamma}{n^d} \sum_{x \in \tdn}  \big( J^n \ast u^k \big)_x w_x^i f \big( \tfrac{x}{n} \big)^2 \Big\} d \mu_t^n \Big)\\
		&\leq \gamma^{-1} \Big( H_n(t) + \log 2 + \frac{1}{8n^{2d}} \sum_{x \in \tdn} \gamma^2 \big( J^n \ast u^k\big)_x^2 f \big( \tfrac{x}{n} \big)^4\Big)\\
		&\leq \gamma^{-1} \Big( H_n(t) + \log 2 + \frac{\gamma^2 \norm{J}_\infty^2 \norm{f}_\infty^4}{n^d} \Big)\\
		&\leq \frac{Ce^{Ce^{Ct}}\norm{f}_\infty^2}{n^{d/2}},
\end{split}
\]
where the last estimate was obtained taking $\gamma = n^{d/2}$.
The third term in \eqref{illapel} can be estimated in the same way. We conclude that
\[
\bb E_{\mu_{\vec u_0}^n}^n \Big[ \Big| \Gamma_n \big( X_s^{n,i}(f), X_s^{n,i+1}(f) \big)
		- \frac{1}{n^d} \sum_{x \in \tdn} \gamma_x^{i,i+1} f \big( \tfrac{x}{n} \big)^2 \Big| \Big] \leq \frac{Ce^{Ce^{Ct}}\norm{f}_\infty^2}{n^{d/2}}.
\]
Bounds of the same form also hold for the other quadratic covariations. 

Up to here, we have proved the lemma with $\gamma^{i,j}(t,x)$ replaced by $\gamma_x^{n,i,j}(t)$. Using Lemma \ref{la6}, we can replace $\gamma_x^{n,i,j}(t)$ by $\gamma^{i,j}(t,x)$, at a cost of order $\frac{Ct e^{Ct}}{n^2}$, which proves the lemma. 
\end{proof}

\subsection{The drift term} 
\label{s4.2}
In this section we compute the drift term
\[
 \int_0^t \big( \partial_{\!s} + L_n \big) X_s^{n,i}(f^i) ds
\]
in the martingale decomposition \eqref{dynkin}. Recall definition \eqref{losvilos}. Observe that $X_t^{n,i}(f)$ can be rewritten as
\[
X_t^{n,i}(f) = \frac{1}{n^{d/2}} \sum_{x \in \tdn} w_x^i f \big( \tfrac{x}{n} \big). 
\]
For $g: \tdn \to \bb R$, we will use the notation
\[
X_t^{n,i}(g) := \frac{1}{n^{d/2}} \sum_{x \in \tdn} w_x^i g_x.
\]
Observe that this notation is coherent with the definition of $\vec X_t^n$ if we identify the function $f \in \mc C^\infty(\bb T^d)$ with its projection $x \mapsto f ( \frac{x}{n})$. In order to further simplify the notation, we will write $f_x$ instead of $f(\frac{x}{n})$.

We have that
\[
L_n \sigma_x^i = \sum_{j \in Q_k} \Big( A_{ij}  + \big( J^n \ast \sigma^k \big)_x M_{ij} \Big) \sigma_x^j.
\]
Therefore,
\[
\begin{split}
\big( \partialt + L_n \big) (\sigma_x^i - u_x^i) 
		&= \sum_{j \in Q_k} \Big( A_{ij}  + \big( J^n \ast u^k \big)_x M_{ij} \Big) w_x^j
			+ \sum_{j \in Q_k} \big( J^n \ast w^k \big)_x M_{ij} u_x^j\\
		&+\sum_{j \in Q_k}  \big( J^n \ast w^k \big)_x M_{ij} w_x^j.
\end{split}
\]
The idea is that the first two terms on the right-hand side of this relation give rise to the two drift terms in \eqref{SPDE}, while the third term gives rise to an error term. For $g: \tdn \to \bb R$, define $J^{n,\ast} \ast g: \tdn \to \bb R$ as
\[
\big( J^{n,\ast} \ast g \big)_x := \frac{1}{n^d} \sum_{y \in \tdn} J_{y,x}^n g_y
\]
for every $x \in \tdn$.
From the computation above,
\begin{equation}
\label{quintero}
\begin{split}
\big( \partialt + L_n \big) X_t^{n,i}(f)  
		&= \sum_{j \in Q_k} \Big\{ A_{ij} X_t^{n,j}(f) 
		+ M_{ij} X_t^{n,j} \big( ( J^n \ast u^k ) f \big)
		+X_t^{n,k}\big( J^{n,\ast} \ast (M_{ij} u^j f)\big)
\\
		&\quad + \frac{1}{n^{3d/2}} \sumxy J_{x,y}^n w_y^k M_{ij} w_x^j f\big( \tfrac{x}{n}\big)\Big\}.
\end{split}
\end{equation}
The first line on this expression corresponds to the discretization of the drift term in \eqref{SPDE}, while the second term is an error term. Let us show that the error term converges to $0$ as $n \to \infty$. By \eqref{donvar2}, for every $\gamma >0$,
\[
\begin{split}
\bb E_{\mu_{\vec u_0}^n}^n \Big[ \Big|  \frac{1}{n^{3d/2}} \sumxy J_{x,y}^n w_y^k w_x^j f\big( \tfrac{x}{n}\big)\Big|\Big] 
		&\leq \gamma^{-1} \Big( H_n(t) + \log 2 +\\
		&\quad + \max_{\pm} \; \log \int \exp\Big\{  \frac{\pm \gamma}{n^{3d/2}} \sumxy J_{x,y}^n w_y^k w_x^j f\big( \tfrac{x}{n}\big)\Big\} d \mu_t^n\Big).
\end{split}
\]
By Hanson-Wright inequality (see Proposition \ref{pa2}), for every 
\[
\gamma^{-1} \geq \Big( \frac{64 \norm{J}_\infty^2 \norm{f}_\infty^2 }{n^{d}}\Big)^{1/2},
\]
\[
\bb E_{\mu_{\vec u_0}^n}^n \Big[ \Big|  \frac{1}{n^{3d/2}} \sumxy J_{x,y}^n w_y^k w_x^j f\big( \tfrac{x}{n}\big)\Big|\Big] 
		\leq \gamma^{-1} \big( H_n(t) + \log 6 \big)
\]
from where
\[
\bb E_{\mu_{\vec u_0}^n}^n \Big[ \Big|  \frac{1}{n^{3d/2}} \sumxy J_{x,y}^n w_y^k w_x^j f\big( \tfrac{x}{n}\big)\Big|\Big] 
		\leq \frac{Ce^{Ce^{Ct}}\norm{f}_\infty}{n^{d/2}},
\]
which proves that this expression is an error term.

By Remark \ref{r1}, we see that
\[
\bb E_{\mu_{\vec u_0}^n}^n \big[ X_t^{n,i}(g)^2 \big] \leq Ce^{Ce^{Ct}} \| g\|_{\ell^2_n}^2.
\]
Using Lemma \ref{la6} together with this estimate and the computations above, we have proved the following result:

\begin{lemma}
\label{l3} For every $f \in \mc C^\infty(\bb T^d; \bb R^{Q_k})$ and every $T >0$,
\[
\lim_{n \to \infty} \int_0^t \Big\{\big( \partial_{\!s} + L_n \big) \vec X_s^n(\vec f) - \vec X_s^n( A^\ast \vec f) 
		- \vec X_s^n\big( (J \ast u^k(s))M^\ast \vec f \big) 
		- X_s^{n,k} \big(J^\ast \ast \< \vec f, M \vec u_s\> \big) \Big\} ds =0
\]
in $\mc L^1(\bb P_{\mu_{\vec u_0}^n}^n)$.
\end{lemma}

\subsection{Tightness}

In this section our aim is to prove tightness of the sequence of processes $(\vec X_t^n; t \geq 0)_{n \in \bb N}$ with respect to the $J_1$-Skorohod topology of \cadlag paths with values in the space of distributions. By Mitoma's criterion \cite{Mit}, it is enough to show tightness of the sequence $(X_t^{n,i}(f); t \geq 0)_{n \in \bb N}$ for every $f \in \mc C^\infty(\bb T^d)$ and every $i \in Q_k$. By decomposition \eqref{dynkin}, it is enough to prove tightness of each of the sequences
\[
\big(X_0^{n,i}(f) \big)_{n \in \bb N}, \big( M_t^{n,i}(f); t \geq 0\big)_{n \in \bb N}
\]
and
\[
\Big(\int_0^t (\partial_{\!s} + L_n) X_s^{n,i}(f) ds; t \geq 0 \Big)_{n \in \bb N}.
\]
The initial condition is the simplest; we have already mentioned that $\vec X_0^n$ converges to a Gaussian process $\vec \xi$, and in particular it is tight. The next in simplicity is the integral term.
Let us use the notation
\[
\vec w_x := \big(w_x^i; i \in Q_k \big).
\]
With respect to the measure $\mu_t^n$, the random variables $(\vec w_x ; x \in \tdn)$ are independent and bounded. 
Therefore, using \eqref{quintero} and McDiarmid's inequality, we see that there exists a finite constant $C=C(a)$ such that
\[
\log \int \exp\big\{ \gamma (\partialt + L_n) X_t^{n,i}(f) \big\} d \mu_t \leq C \norm{f}_\infty (1+ \norm{J}_\infty) \gamma^2.
\]
for every $\gamma \in \bb R$. Therefore,
\[
\log \int \exp\big\{ \gamma \big( (\partialt + L_n) X_t^{n,i}(f) \big)^2 \big\} d \mu_t \leq \log 3
\]
for\[
\gamma^{-1} = 8 C^2 \norm{f}_\infty^2 (1+ \norm{J}_\infty)^2.
\]

Observe that in Section \ref{s4.2} we showed that this term is well approximated in $\mc L^1(\bb P_{\mu_{\vec u_0}^n}^n)$ by a linear function of the fluctuations. However, the estimates of Section \ref{s4.2} do not yield good bounds in $\mc L^2(\bb P_{\mu_{\vec u_0}^n}^n)$. For this reason, here we need to do something different.

By \eqref{donvar}, for every $\gamma >0$, 
\[
\bb E_{\mu_{\vec u_0}^n}^n \Big[ \Big( \big( \partialt + L_n\big) X_t^{n,i}(f) \Big)^2 \Big]
		\leq \gamma^{-1} \Big( H_n(t) +
\log \int \exp\big\{ \gamma \big( (\partialt + L_n) X_t^{n,i}(f) \big)^2 \big\} d \mu_t \Big),
\]
from where we obtain the bound
\begin{equation}
\label{valparaiso}
\bb E_{\mu_{\vec u_0}^n}^n \Big[ \Big( \big( \partialt + L_n\big) X_t^{n,i}(f) \Big)^2 \Big]  \leq Ce^{Ce^{Ct}} \norm{f}^2_\infty(1+ \norm{J}_\infty)^2.
\end{equation}
By Cauchy-Schwartz inequality, we see that there exists a constant $C = C(a,T,\norm{J}_\infty, \norm{f}_\infty)$ such that for every $t_1,t_2 \in [0,T]$ such that $t_1 < t_2$,
\[
\bb E_{\mu_{\vec u_0}^n}^n \Big[ \Big( \int_{t_1}^{t_2} (\partialt + L_n) X_t^{n,i}(f) dt \Big)^2 \Big] \leq C(t_2-t_1)^2.
\]
We conclude by Kolmogorov-Centsov's criterion that the sequence 
\[
\Big(\int_0^t (\partial_{\!s} + L_n)  X_s^{n,i}(f) ds; t \geq 0 \Big)_{n \in \bb N}
\]
is tight with respect to the uniform topology on continuous paths, as we wanted to show.

Now we are left to prove tightness of the martingale sequence 
\[
\big( M_t^{n,i}(f); t \geq 0 \big)_{n \in \bb N}.
\]
We will use the following criterion (see Theorem 3.6 in \cite{Whi} and also Theorem VI.4.13 in \cite{JacShi}):

\begin{proposition}
\label{p2} The sequence $\big((M_t^{n,i}(f); i \in Q_k); t \geq 0 \big)_{n \in \bb N}$ is tight if for every $i \in Q_k$, the sequence $\big( \<M^i(f), M^i(f)\>_t; t \geq 0 \big)_{n \in \bb N}$ is $\mc C$-tight.
\end{proposition}

From \eqref{vicuna}, we see that 
\[
\Gamma_n \big( X_t^{n,i}(f), X_t^{n,i}(f) \big) \leq  \big( a + \norm{J}_\infty \big) \norm{f}_{\ell_n^2}^2.
\]
Therefore, $(\<M^i(f), M^i(f)\>_t; t \geq 0)$ is Lipschitz, uniformly in $t$ and $n$, and in particular it is $\mc C$-tight, which proves tightness of the martingale processes, and in consequence of the processes $(\vec X_t^n; t \geq 0)$.

\subsection{Convergence}

In this section we will prove Theorem \ref{t2}. The proof follows the usual strategy of convergence of stochastic processes. In the previous section, we have proved tightness of $(\vec X_t^n; t \geq 0)_{n \in \bb N}$ with respect to the $J_1$-Skorohod topology of \cadlag paths in the space of distributions. Therefore, convergence follows if we prove uniqueness of accumulation points of this sequence. Let $n'$ be a subsequence of $n$ for which $(\vec X_t^{n'}; t \geq 0)_{n'}$ converges weakly to a process $(\vec X_t; t \geq 0)$. First we observe that this convergence and Lemma \ref{l3} imply that
\[
\lim_{n \to \infty} \int_0^t (\partial_{\!s} + L_n) \vec X_s^{n'}(\vec f) ds
		= \int_0^t \Big( \vec X_s( A^\ast \vec f) + \vec X_s( ( J \ast u_s^k ) M^\ast f) + X_s^k\big( J^\ast \ast \<\vec f, M \vec u_s\>\big) \Big) ds
\]
for every $t \geq 0$ and every $f\in \mc C^\infty(\bb T^d; \bb R^{Q_k})$. Observe however, that due to the structure of the $J_1$-Skorohod topology, the convergence of $(\vec X_t^{n'}; t \geq 0)$ \emph{does not imply} the convergence of $\vec X_t^{n'}(\vec f)$ for a fixed time $t$. It does imply, nevertheless, the convergence of $\vec X_t^{n'}(\vec f)$ to $\vec X_t(\vec f)$ for \emph{almost every} time $t$. To be more precise, weak convergence with respect to the $J_1$-Skorohod topology implies convergence of fixed-times marginals except for a set of times which is at most countable. The condition under which the convergence holds for all times $t \geq 0$, is \emph{stochastic continuity} of the limiting process. Let us prove that $(\vec X_t; t \geq 0)$ has continuous paths, so in particular it is stochastically continuous.

Let us define the process $(\Delta_t^n(f); t \geq 0)$ as
\[
\Delta_t^n(f) := \sup_{0 \leq s \leq t} \| \vec M_s^n(f) - \vec M_{s-}^n(f)\|
\]
for every $t \geq 0$.
In other words, $\Delta_t^n(f)$ is the size of the largest jump of the martingale $(\vec M_t^n(f); t \geq 0)$ up to time $t$.
We have the deterministic bound
\begin{equation}
\label{quillota}
\Delta_t^n(f) \leq \frac{2 \norm{f}_\infty}{n^{d/2}}.
\end{equation}
By upper semicontinuity, we conclude that every limit point of $(\vec M_t^n; t \geq 0)_{n \in \bb N}$ is continuous. The integral process in \eqref{dynkin} is continuous to start with, so in particular every limit point is also continuous, from where we conclude that every limit point of $(\vec X_t^n; t \geq 0)$ is continuous. In particular, we have proved that
\[
\lim_{n' \to \infty} \vec X_t^{n'} = \vec X_t \text{ for every } t \geq 0.
\]
Let us show the convergence of the martingale processes $(\vec M_t^n; t \geq 0)_{n \in \bb N}$. We will see that we do not need to take subsequences in this case. First we observe that by Wald's device, convergence of $(\vec M_t^n; t \geq 0)_{n \in \bb N}$ to a process $(\vec M_t; t \geq 0)$ is equivalent to convergence of $(\vec M_t^n(\vec f); t \geq 0)_{n \in \bb N}$ to $(\vec M_t(\vec f); t \geq 0)$ for every $\vec f \in \mc C^\infty(\bb T^d; \bb R^{Q_k})$.  By Theorem VIII.3.11 of \cite{JacShi} and \eqref{quillota}, convergence of $(\vec M_t^n(\vec f); t \geq 0)$ to a Gaussian martingale $(\vec M_t( \vec f); t \geq 0)$ is equivalent to convergence of the quadratic covariations $(\< M^{n,i}(f), M^{n,j}(g)\>_t; t \geq 0)_{n \in \bb N}$ to the corresponding limits $(\<M^i(f),M^j(g)\>_t; t \geq 0)$. But this is exactly what Lemma \ref{l2} proves. Therefore, we have proved the following result:

\begin{lemma}
\label{l4}
The sequence $(\vec M_t^n; t \geq 0)_{n \in \bb N}$ converges in distribution with respect to the $J_1$-topology to the unique Gaussian martingale process $(\vec M_t; t \geq 0)$ satisfying $\vec M_0 = \vec 0$ and
\[
\<M^i(f),M^j(f)\>_t =  \int_0^t \int \gamma^{i,j}(s,x) f(x)^2 dx ds
\]
for every $f \in \mc C^\infty(\bb R^d)$, every $i,j \in Q_k$ and every $t \geq 0$.
\end{lemma}

Observe that $\vec X_0^n$ converges to $\vec \xi$. Taking limits along $n'$ in \eqref{dynkin}, we have proved that
\begin{equation}
\label{rancagua}
\vec X_t(\vec f) = \vec \xi(\vec f) + \int_0^t \Big(  \vec X_s(A^\ast \vec f)  + \vec X_s( ( J \ast u_s^k ) M^\ast \vec f) + X_s^k\big( J^\ast \ast \< \vec f, M \vec u_s\>\big) \Big) ds
		+ \vec M_t(\vec f)
\end{equation}
for every $t \geq 0$ and every $\vec f \in \mc C^\infty(\bb T^d; \bb R^{Q_k})$. 
From \eqref{valparaiso}, we also see that
\begin{equation}
\label{rengo}
\bb E\big[ \norm{\vec X_t(\vec f)} \big] \leq Ce^{Ce^{Ct}} \norm{\vec f}_\infty(1+ \norm{J}_\infty)
\end{equation}
for every $t \geq 0$ and every $\vec f \in \mc C^\infty(\bb T^d; \bb R^{Q_k})$.
We claim that \eqref{rancagua} and \eqref{rengo} characterize the law of $(\vec X_t ; t \geq 0)$. First we observe that using \eqref{rengo} we can take suitable approximations to show that for every $t \geq 0$ and every $f \in \mc C([0,t]; \mc C^\infty(\bb T^d; \bb R^{Q_k}))$,
\begin{equation}
\label{marchihue}
\begin{split}
\vec X_t(\vec f_t) 
		= \vec \xi(\vec f_0)
		&+ \int_0^t \Big( \vec X_s(\partial_{\!s} \vec f_s) +  \vec X_s(A^\ast \vec f_s) + \vec X_s( ( J \ast u_s^k ) M^\ast \vec f_s)\\
		&\quad \quad \quad  + X_s^k\big( J^\ast \ast \<\vec f_s, M \vec u_s\> \Big) ds
		+ \int_0^t d \vec M_s(\vec f_s).
\end{split}
\end{equation}
Given a function $\vec f \in \mc C^\infty(\bb T^d; \bb R^{Q_k})$ and $t \geq 0$, let $(P_s \vec f; s \in [0,t])$ be the solution of the backwards Fokker-Planck equation
\begin{equation}
\label{fokpla}
\partial_{\!s} \vec g_s + A^\ast \vec g_s + \big( J \ast u^k_s \big) M^\ast \vec g_s +  \big( J^\ast \ast \<\vec g_s, M \vec u_s \>\big) e_k =0
\end{equation}
with final condition $\vec g_t = \vec f$. Observe that $(P_{t-s} \vec f; s \in [0,t])$ satisfies a parabolic equation of the form  \eqref{curico}. Therefore, Lemma \ref{la5} also applies for \eqref{fokpla}, and in particular $(P_s \vec f; s \in [0,t])$ exists and it is unique.
Using $(P_s \vec f; s \in [0,t])$ as a test function in \eqref{marchihue}, we see that
\begin{equation}
\label{mild}
\vec X_t(\vec f) = \vec \xi(P_0 \vec f) + \int_0^t d \vec M_s( P_s \vec f), 
\end{equation}
which characterize the one-dimensional marginals of $(\vec X_t; t \geq 0)$. Let $s,t \geq 0$ be such that $s<t$. Observe that
\[
\vec X_t(\vec f) = \vec X_s(P_s \vec f) + \int_s^t d \vec M_u( P_u \vec f).
\]
The law of $\vec X_s$ has already been characterized, and $\int_s^t d \vec M_u( P_u \vec f)$ is independent of $\vec X_s$. Therefore, the two-dimensional marginals of $(\vec X_t; t \geq 0)$ are also characterized by this relation.
Recursively, we see that \eqref{mild} characterizes all finite-dimensional marginals of $(\vec X_t; t \geq 0)$. We conclude that the only accumulation point of $(\vec X_t^n; t \geq 0)_{n \in \bb N}$ is the unique solution of \eqref{SPDE}, from where convergence follows.

\appendix

\section{The solution theory of \eqref{echid}}

\label{a1}

\subsection{Wellposedness of \eqref{echid}}

\label{sa1}

In this section we discuss the existence and uniqueness of solutions of the hydrodynamic equation \eqref{echid}, as well as their regularity. It turns out that \eqref{echid} can be understood as an \emph{evolution equation}, that is, as an ODE in a functional space. Since \eqref{echid} has a quadratic nonlinearity, some care is needed in order to solve it. 

As usual in the theory of PDEs, we treat \eqref{echid} as a fixed-point problem for an inhomogeneous linear equation. For $\vec u: \bb T^d \to \bb R^{Q_k}$ and $v: \bb T^d \to \bb R$, let $\Lambda(\vec u, v) : \bb T^d \to \bb R^{Q_k}$ be given by
\[
\Lambda(\vec u, v)(x) := A \vec u(x) + \big(J \ast v (x) \big) M \vec u(x) 
\]
for every $x \in \bb T^d$. Take $T>0$ and $\vec u_0: \bb T^d \to \bb R^{Q_k}$. Given $v: \bb T^d \times [0,T] \to \bb R$ bounded, we say that $\vec u^v=( \vec u^v(t,x); t \in [0,T], x \in \bb T^d)$ is a solution of
\begin{equation}
\label{echid2}
\partialt \vec u_t = \Lambda(\vec u_t, v_t)
\end{equation}
with initial condition $\vec u_0$ if
\begin{equation}
\label{intode}
\vec u^v_t = \vec u_0 + \int_0^t \Lambda(\vec u_s^v; v_s) ds
\end{equation}
for every $t \in [0,T]$.
Observe that this equation is linear in $\vec u$, and therefore boundedness of $v$ implies existence and uniqueness of solutions. 

Let $\norm{A}$, $\norm{M}$ be the operator norms of the matrices $A$, $M$, that is,
\[
\norm{A} := \sup_{\norm{\vec u}=1} \norm{ A \vec u}, \quad \norm{M} := \sup_{\norm{\vec u}=1} \norm{ M \vec u}.
\]
Define as well
\[
\norm{J}_1 := \sup_{x \in \bb T^d} \int J(x,y) dy
\]
and
\[
\lambda_t(v) := \sup_{0 \leq s \leq t} \big( \norm{A} + \norm{J}_1 \norm{M} \norm{v_s}_\infty\big).
\]
By Holder's inequality,
\[
\norm{ J \ast v_t}_\infty \leq \norm{J}_1 \norm{v_t}_\infty
\]
We have the following estimate:

\begin{lemma}
\label{la1}
For every $v: \bb T^d \times [0,T]$ bounded and every $t \in [0,T]$, 
\[
\norm{\vec u^v_t}_\infty \leq \norm{\vec u_0}_\infty e^{\lambda_t(v) t}.
\]
\end{lemma}

\begin{proof}
Observe that
\[
\norm{\Lambda(\vec u, v_t)} \leq \big( \norm{A} + \norm{J}_1 \norm{v_t}_\infty \norm{M} \big) \norm{\vec u} \leq \lambda_t(v) \norm{\vec u}.
\]
Therefore, for every $s \in [0,t]$,
\[
\norm{\vec u_s^v} \leq \norm{\vec u_0} +\int_0^s \lambda_{s'}(v) \norm{\vec u_{s'}^v} ds'
		\leq \norm{\vec u_0} +\lambda_t(v) \int_0^s  \norm{\vec u_{s'}^v} ds'
\]
and by Gronwald's inequality,
\[
\norm{\vec u_t^v} \leq \norm{\vec u_0} e^{\lambda_t(v)t},
\]
as we wanted to show.
\end{proof}

This lemma serves as an \emph{a priori} bound for solutions of \eqref{echid2}. The following lemma states the continuity with respect to parameters of solutions of \eqref{echid2}:

\begin{lemma}
\label{la2} 
For $v$, $w$ bounded,
\[
\norm{\vec u_t^v - \vec u_t^w} \leq \norm{J}_1 \norm{M} \norm{\vec u_0}_\infty e^{(\lambda_t(v)+\lambda_t(w))t} \int_0^t \norm{v_s-w_s}_\infty ds
\]
for every $t \in [0,T]$.
\end{lemma}

\begin{proof}
Observe that for every $\vec u: \bb T^d \to \bb R^{Q_k}$ and every $t \in [0,T]$,
\[
\norm{\Lambda(\vec u, v_t) - \Lambda(\vec u, w_t)}_\infty \leq \norm{J}_1 \norm{v_t-w_t}_\infty \norm{M} \norm{\vec u}_\infty.
\]
Therefore,
\[
\begin{split}
\norm{\Lambda(\vec u_t^v, v_t) - \Lambda(\vec u_t^w, w_t)}_\infty \leq \lambda_t(v) \norm{\vec u_t^v - \vec u_t^w}_\infty +  \norm{J}_1 \norm{v_t-w_t}_\infty \norm{M} \norm{\vec u_t^w}_\infty.
\end{split}
\]
Inserting this estimate into \eqref{intode}, we see that
\[
\begin{split}
\norm{\vec u_t^v - \vec u_t^w}_\infty 
		&\leq \int_0^t \norm{\Lambda(\vec u_s^v, v_s) - \Lambda(\vec u_s^w,w_s)}_\infty ds \\
		&\leq \lambda_t(v) \int_0^t \norm{\vec u_s^v - \vec u_s^w}_\infty ds 
				+\norm{J}_1 \norm{M} \norm{\vec u_t^w}_\infty  \int_0^t \norm{v_s-w_s}_\infty ds\\
		&\leq \lambda_t(v) \int_0^t \norm{\vec u_s^v - \vec u_s^w}_\infty ds 
				+\norm{J}_1 \norm{M} \norm{\vec u_0}_\infty e^{\lambda_t(w)t} \int_0^t \norm{v_s-w_s}_\infty ds.
\end{split}
\]
By Gronwald's inequality, we conclude that 
\[
\norm{\vec u_t^v - \vec u_t^w}_\infty \leq \norm{J}_1 \norm{M} \norm{\vec u_0}_\infty e^{(\lambda_t(v) + \lambda_t(w))t} \int_0^t \norm{v_s- w_s}_\infty ds,
\]
which proves the lemma.
\end{proof}

Observe that \eqref{echid} can be written as
\[
\partialt \vec u_t = \Lambda( \vec u_t, u^k).
\]
Let $\vec u_0: \bb T^d \to \mc P_{k}$ be fixed. For each bounded function $\vec v: \bb T^d \to \bb R^{Q_k}$ such that $\vec v_0 = \vec u_0$, let $\mc R \vec v: \bb T^d \to \bb R^{Q_k}$ be the solution of \eqref{echid2} with initial condition $\vec u_0$. In other words, $\mc R \vec v = \vec u^{v^k}$. Equation \eqref{echid} can be written as a fixed-point problem:
\[
\vec u = \mc R \vec u.
\]
Our aim is to find a suitable space on which the operator $\mc R$ is a contraction, on which case this fixed-point equation has a unique solution. We have the following result:

\begin{lemma}
\label{la3}
There exists $T_0 = T_0(\norm{A}, \norm{M}, \norm{J}_1)>0$ such that the equation 
\[
\partialt \vec u_t = \Lambda(\vec u_t, u^k_t)
\]
has a unique solution in $\bb T^d \times [0,T_0]$ for every initial condition $\vec u_0: \bb T^d \to \mc P_{k}$.
\end{lemma}

\begin{proof}
Observe that $\norm{\vec u_0}_\infty \leq 1$. For $\alpha, T>0$, define
\[
V_{\alpha,T}:= \big\{\vec v: \bb T^d \times [0,T] \to \bb R^{Q_k}; \lambda_T(v) \leq \alpha \big\}
\]
and denote by $\tnorm{\cdot}$ in $V_{\alpha,T}$ the supremum norm. 
By Lemma \ref{la1},
\[
\begin{split}
\lambda_T\big( \big(\mc R \vec v\big)^k\big) 
		&\leq \sup_{0 \leq t \leq T} \big( \norm{A} + \norm{J}_1 \norm{M} \norm{(\mc R \vec v)_{\!t}}_\infty \big) \\
		&\leq  \norm{A} + \norm{J}_1 \norm{M} \norm{\vec u_0}_\infty e^{\lambda_T(v^k)T}.
\end{split}
\]
Therefore, given $\alpha > \norm{A}$, we see that $\lambda_T((\mc R \vec v)^k) \leq \alpha$ whenever $\lambda_T(v^k) \leq \alpha$ and
\[
T \leq \frac{1}{\alpha} \log \Big( \frac{\alpha - \norm{A}}{\norm{J}_1 \norm{M}}\Big) =: T(\alpha).
\]
In particular, $\mc R$ maps $V_{\alpha,T}$ into $V_{\alpha,T}$ for $T \leq T(\alpha)$. 

Let us check now under which conditions we can guarantee that $\mc R$ is a contraction. Lemma \ref{la2} tells us that
\[
\norm{(\mc R \vec v)_{\!t} - (\mc R \vec w)_{\!t}}_\infty 
		\leq \norm{J}_1 \norm{M} e^{2 \alpha t} \int_0^t \norm{v_s^k - w_s^k}_\infty ds.
\]
We conclude that
\[
\tnorm{ \mc R \vec v -\mc R \vec w} \leq  \norm{J}_1 \norm{M} T e^{2 \alpha T} \tnorm{\vec v- \vec w}.
\]
The multiplicative factor  $\norm{J}_1 \norm{M} T e^{2 \alpha T}$ converges to $0$ as $T \to 0$, and therefore it is smaller than $1$ for $T$ small. This proves the lemma. 
\end{proof}

\begin{remark}
In the previous lemma, it is possible to find a more concrete value for $T_0$. In fact, assuming that $T \leq T(\alpha)$, we see that
\[
 \norm{J}_1 \norm{M} T e^{2 \alpha T} \leq  \norm{J}_1 \norm{M}T e^{2 \alpha T(\alpha)}
 		= \frac{T(\alpha-\norm{A})^2}{\norm{J}_1 \norm{M}},
\]
which is less or equal than $1/2$ for
\[
T = \frac{\norm{J}_1 \norm{M}}{2(\alpha-\norm{A})^2}.
\]
Therefore, it is enough to take 
\[
T_0 = \min\Big\{ \frac{\norm{J}_1 \norm{M}}{2(\alpha-\norm{A})^2}, \frac{1}{\alpha} \log \Big( \frac{\alpha - \norm{A}}{\norm{J}_1 \norm{M}}\Big)\Big\}.
\]
\end{remark}

Lemma \ref{la3} only gives \emph{local} existence and uniqueness of solutions for \eqref{echid}. However, for the type of solutions of \eqref{echid} we are interested, one can actually show \emph{global} existence and uniqueness, as the following lemma shows:

\begin{lemma}
\label{la4}
Let $\vec u_0: \bb T^d \to \mc P_k$ be such that $u_0^i \geq \epsilon_0 >0$ for every $i \in Q_k$. For every $t \in [0,T_0]$ and every $i \in Q_k$,
\[
\vec u_t^i(x) >0.
\] 
In particular, the solution of \eqref{echid} can be extended to $t \in [0,\infty)$ in a unique way.
\end{lemma}

\begin{proof}
Let $\partial \mc P_k$ denote the boundary of $\mc P_k$. Observe that the lemma can be rephrased as
\[
\text{If } \vec u_0 \notin \partial \mc P_k, \text{ then } \vec u_t \notin \partial \mc P_k \text{ for every } t \in [0,T_0].
\]
Define
\begin{equation}
\label{ancud}
T_1:= \inf \big\{t \geq 0; \vec u_t(x)  \in \partial \mc P_k  \text{ for some } x \in \bb T^d\big\},
\end{equation}
where we take $T_1 := T_0$ if there is no such time. Up to $T_1$, all coordinates of $\vec u_t$ are positive, and $u^k_t$ is smaller than $1$. Therefore,
\[
\partialt u_t^i(x) \geq 
\left\{
\begin{array}{c@{\;;\;}l}
- a u_t^i(x) & i=k, \\
- \norm{J}_1 u_t^i(x) & i\neq k.
\end{array}
\right.
\]
We conclude that
\begin{equation}
\label{castro}
u_t^i(x) \geq \epsilon_0 e^{- \max\{a,\norm{J}_1\}t}
\end{equation}
for every $i \in Q_k$, every $x \in \bb T^d$ and every $t \leq T_1$. By continuity, we see that the infimum in \eqref{ancud} runs over an empty set and therefore $T_1=T_0$. Moreover, since $\sum_i u_t^i=1$ for every $t \in [0,T_0]$, we conclude that $\vec u_{T_0} \notin \partial \mc P_k$. In consequence, we can use Lemma \ref{la3} to extend the definition of $\vec u_t$ to the interval $[0,2T_0]$, and use the argument above to argue that $T_1 = 2 T_0$. Recursively, we can extend the solution of \eqref{echid} to $t \in [0,\infty)$, for which $T_1 =+\infty$, as we wanted to show.
\end{proof}

\subsection{Conservation of regularity}

In this section we show that smooth initial conditions with values in $\mc P_k$ give rise to smooth solutions of \eqref{echid}. We point out that \eqref{echid} does not have an obvious regularizing effect, so the best we can do without a deeper understanding of \eqref{echid}, is to prove that smooth solutions remain smooth. We will prove the following result: 

\begin{lemma} 
\label{la5}
Let $\ell \in \bb N$ and let $\vec u_0: \bb T^d \to \mc P_k$ be of class $\mc C^\ell$. Define $\alpha := \norm{A} + 2 \norm{J}_1 \norm{M}$. There exists a finite constant 
\[
C_\ell = C_\ell(\norm{A} \norm{J}_1, \norm{M}, \norm{\vec u_0}_{\mc C^\ell})
\]
such that
\[
\norm{\partial_{\vec{ I}}^\ell \vec u_t}_\infty \leq C_\ell e^{\ell \alpha t} 
\]
for every multiindex $\vec I$ of order $\ell$ and every $t \geq 0$.
\end{lemma}

\begin{proof}
We proceed by induction on $\ell$. By Lemma \ref{la4}, $\norm{\vec u_t}_\infty \leq 1$, and the lemma holds for $\ell =0$. Let $\ell \geq 1$, assume that the lemma holds for $\ell' < \ell$ and let $\vec u_0$ be of class $\mc C^\ell$. Let $\vec I$ be a multiindex of order $\ell$ and let $\vec g_t := \partial^\ell_{\vec I} \vec u_t$. If $\vec g_t$ is well defined on an interval of the form $[0,T_\ell]$, then it satisfies the equation
\begin{equation}
\label{curico}
\begin{split}
\partialt \vec g_t 
		&= A \vec g_t + \big( J \ast u_t^k \big) M \vec g + \big( J \ast g_t^k \big) M \vec u_t\\
		&\quad + \sum_{i=1}^{\ell-1} \sum_{\vec I_1, \vec I_2} \big( J \ast \partial^{i}_{\vec I_1} u_t^k \big) M \big(\partial^{\ell-i}_{\vec I_2} \vec u_t \big),
\end{split}
\end{equation}
where the sum runs over all the $\binom{\ell}{i}$ ways to split the multiindex $\vec I$ into two multiindices $\vec I_1$, $\vec I_2$ of orders $i$, $\ell-i$ respectively. Observe that this equation is linear in $\vec g_t$, and that by the inductive hypothesis, the coefficients of this equation are bounded. Therefore this equation has a unique solution, and it satisfies
\[
\begin{split}
\norm{\vec g_t}_\infty 
		&\leq \norm{\vec g_0}_\infty + \int_0^t \Big( \big( \norm{A} + 2 \norm{J}_1 \norm{M} \big) \norm{\vec g_s}_\infty + \sum_{i=1}^{\ell-1} \binom{\ell}{i} \norm{J}_1 \norm{M} C_i C_{\ell-i} e^{i \alpha s} e^{(\ell-i) \alpha s}\Big)ds\\
		&\leq \norm{\vec g_0}_\infty + \int_0^t \alpha \norm{\vec g_s}_\infty ds + C_{\ell} \big( e^{\ell \alpha t} -1 \big).
\end{split}
\]
Using the integrating factor $e^{-\alpha t}$, we see that
\[
\norm{\vec g_t}_\infty \leq \norm{\vec g_0}_\infty e^{\alpha t} + \frac{C_\ell}{(\ell-1) \alpha} \big( e^{\ell \alpha t}- e^{\alpha t}\big),
\]
which proves the lemma.
\end{proof}

\subsection{The numerical approximation}
\label{ra1}
In this section we prove that for smooth initial conditions $\vec u_0$ and smooth kernels $J$, the solutions of \eqref{numaprox} converges uniformly to the solution of \eqref{echid}, with an explicit speed of convergence.

It will be convenient to introduce a discretization of the operator $\Lambda$. For $\vec u: \tdn \to \bb R^{Q_k}$ and $v: \tdn \to \bb R$, let $\Lambda^n(\vec u, v): \tdn \to \bb R^{Q_k}$ be given by
\[
\Lambda^n(\vec u, v)_x := A \vec u_x + \big( J^n \ast v\big)\vphantom{v}_x M \vec u_x
\]
for every $x \in \tdn$. Define
\[
\norm{J}_{1,n} := \sup_{x \in \tdn} \frac{1}{n^d} \sum_{y \in \tdn} J_{x,y}^n.
\]
Observe that $\Lambda^n$ satisfies the following estimates:
\begin{equation}
\label{talca}
\norm{\Lambda^n(\vec u, v)_x - \Lambda^n(\vec u, w)_x} \leq \norm{J}_{1,n} \norm{M} \norm{\vec u_x} \norm{v-w}_\infty,
\end{equation}
\begin{equation}
\label{linares}
\norm{\Lambda^n(\vec u^1,v)_x - \Lambda^n(\vec u^2,v)_x} \leq \big( \norm{A} + \norm{J}_{1,n} \norm{M} \norm{v}_\infty \big) \norm{\vec u^1 - \vec u^2}_\infty.
\end{equation}
With these estimates in hand, we can repeat the proofs of Lemmas \ref{la1}-\ref{la4} to obtain the corresponding bounds for solutions of \eqref{numaprox}, replacing $\norm{J}_1$ by $\norm{J}_{1,n}$. We will use these bounds without further ado, and we leave the details to the interested reader.

Let $\vec u: \bb T^d \to \mc P_k$ be smooth. Let $(\vec u(t,x) ; t \geq 0, x \in \bb T^d)$ be the solution of \eqref{echid} with initial condition $\vec u_0$ and let $(\vec u_x(t); t \geq 0, x \in \tdn)$ be the solution of \eqref{numaprox} with initial condition $\vec u_x^n := \vec u_0(\frac{x}{n})$. We have the following result:

\begin{lemma}
\label{la6}
There exists a finite constant $C = C(d,\norm{J}_{\mc C^2},\norm{M})$ such that
\[
\norm{\vec u_t - \vec u_t^n}_\infty \leq \frac{Cte^{Ct}}{n^2}
\]
for every $t \geq 0$ and every initial condition $\vec u_0: \bb T^d \to \mc P_k$.
\end{lemma}

\begin{proof}
The first step is to prove that $\Lambda^n$ is a numerical approximation of $\Lambda$. Let $\vec u: \bb T^d \to \mc P_k$ be a continuous function. We will also denote by $\vec u$ the function $x \mapsto \vec u(\frac{x}{n})$ defined in $\tdn$. We have that
\[
\Lambda(\vec u, u^k)\big(\tfrac{x}{n}\big) - \Lambda^n(\vec u, u^k)_x = \big( J \ast u^k \big(\tfrac{x}{n} \big) - (J^n \ast u^k)_x\big) M \vec u \big( \tfrac{x}{n}\big),
\]
from where 
\[
\norm{ \Lambda(\vec u, u^k) - \Lambda^n(\vec u, u^k)}_\infty \leq \norm{M} \sup_{x \in \tdn} \Big| J \ast u^k \big(\tfrac{x}{n} \big) - (J^n \ast u^k)_x \Big|.
\]
Here we have used that $\norm{\vec u}_\infty \leq 1$.
The trapezoidal rule states that for smooth functions $f: \bb R \to \bb R$ of period $1$, 
\[
\Big| \frac{1}{n} \sum_{x=1}^n f \big( \tfrac{x}{n} \big) - \int_0^1 f(x) dx \Big| \leq \frac{1}{12n^2}\norm{f''}_\infty. 
\]
Using this rule once on each coordinate of $\bb T^d$, we see that
\[
\Big| J \ast u^k \big(\tfrac{x}{n} \big) - (J^n \ast u^k)_x \Big| \leq \frac{1}{8 n^2} \sum_{i=1}^d \norm{\partial^2_{ii} ( J \ast u^k)}_\infty \leq \frac{d}{8 n^2} \norm{D^2J}_\infty 
\]
This bound has the following consequence:
\[
\begin{split}
\norm{\Lambda(\vec u_t, u_t^k) - \Lambda^n(\vec u_u^n, u_t^{n,k})}_\infty
		&\leq \norm{\Lambda(\vec u_t, u_t^k) - \Lambda^n(\vec u_t, u_t^k)}_\infty
				+ \norm{\Lambda^n(\vec u_t, u_t^k) - \Lambda^n( \vec u_t^n, u_t^k)}_\infty \\
		& \quad + \norm{\Lambda^n(\vec u_t^n, u_t^k)-\Lambda^n(\vec u_t^n, u_t^{n,k})}_\infty\\
		&\leq \frac{d \norm{D^2J}_\infty}{8n^2} + \big( \norm{A} + \norm{J}_{1,n} \norm{M} \big) \norm{\vec u_t - \vec u_t^n}_\infty\\
		&\quad +\norm{J}_{1,n} \norm{M} \norm{u_t^k-u_t^{n,k}}_\infty.
\end{split}
\]
Using the integral formulations of equations \eqref{echid} and \eqref{numaprox}, we see that
\[
\begin{split}
\norm{\vec u_t - \vec u_t^n}_\infty 
		&\leq \int_0^t \Big( \frac{d \norm{D^2J}_\infty}{8n^2} +\big( \norm{A} + \norm{J}_{1,n} \norm{M} \big) \norm{\vec u_s - \vec u_s^n}_\infty
			\!+\norm{J}_{1,n} \norm{M} \norm{u_s^k-u_s^{n,k}}_\infty \Big)ds\\
		&\leq \big( \norm{A} + 2  \norm{J}_{1,n} \norm{M} \big)  \int_0^t	\norm{\vec u_s - \vec u_s^n}_\infty ds + \frac{d \norm{D^2J}_\infty t}{8n^2}. 
\end{split}			
\]
By Gronwall's inequality, we conclude that
\[
\norm{\vec u_t - \vec u_t^n}_\infty \leq \frac{d \norm{D^2J}_\infty t e^{(\norm{A}+2 \norm{J}_{1,n} \norm{M}) t}}{8n^2}
\]
Observing that $\norm{J}_{1,n} \leq \norm{J}_\infty$, the lemma is proved.
\end{proof}


\section{Concentration bounds}
\label{apb}
In this section we collect various entropy inequalities used throughout this article. We start with Donsker-Varadhan \emph{entropy inequality}:

\begin{proposition}
\label{pa1}
Let $\mu$ be a probability measure and let $f$ be a density with respect to $\mu$. For every function $g$ and every $\gamma>0$,
\begin{equation}
\label{donvar}
\int g f d \mu \leq \gamma^{-1} \Big( \int f \log f d \mu  + \log \int e^{\gamma g} d \mu\Big).
\end{equation}
\end{proposition}

A very simple trick allows to estimate the integral of $|g|$ using this estimate. Observe that for every $x \in \bb R$, $e^{|x|} \leq 2 \max\{e^x, e^{-x}\}$. Combining this estimate with \eqref{donvar}, we obtain the bound
\begin{equation}
\label{donvar2}
\int |g|f d \mu \leq \gamma^{-1} \Big( \int f \log f d\mu + \max_{\pm}\log \int e^{\pm \gamma g} d\mu + \log 2\Big).
\end{equation}

Next we introduce the concept of \emph{sub-Gaussian random variable}. Most of the material on sub-Gaussian random variables exposed here can be found in \cite{Ver}
We say that a real-valued random variable $X$ is sub-Gaussian of index $\sigma^2$ if
\[
\bb E[e^{\theta X}] \leq e^{\frac{1}{2} \sigma^2 \theta^2 }
\]
for every $\theta \in \bb R$. For a given sub-Gaussian random variable $X$, let 
\[
\|X\|_{\psi_2}^2 := \sup_{\theta \in \bb R} \frac{2}{\theta^2} \log \bb E[e^{\theta X}]
\]
be the smallest of such constants. It turns out that $\|\cdot\|_{\psi_2}$ is a norm on the space of sub-Gaussian random variables. 
The squared norm $\|X\|_{\psi_2}^2$ behaves in many ways as the variance of $X$. In particular, we have the following lemmas:

\begin{lemma}
\label{pitagoras}
Let $X_1,...,X_n$ be independent, sub-Gaussian random variables. We have that
\[
\| X_1+ \dots + X_n \|_{\psi_2}^2 \leq \|X_1\|_{\psi_2}^2 + \dots + \|X_n\|_{\psi_2}^2
\]
\end{lemma}

\begin{proof}
Observe that
\[
\begin{split}
\log \bb E\big[ \exp\big\{\theta(X_1+\dots+X_n)\big\}\big] 
		&= \log \bb E[e^{\theta X_1}] + \dots + \log \bb E[e^{\theta X_n}] \\
		&\leq \tfrac{1}{2}\theta^2 \big( \|X_1\|_{\psi_2}^2 + \dots + \|X_n\|_{\psi_2}^2\big),
\end{split}
\]
from where the lemma follows.
\end{proof}

For a proof of the following lemma, see Lemma F.7 in \cite{JarMen}.

\begin{lemma}
\label{quad}
Let $X$ be sub-Gaussian. We have that
\[
\bb E[e^{\gamma X^2}] \leq 3
\]
whenever 
\[
\gamma \leq \frac{1}{4 \|X\|_{\psi_2}^2}.
\]
\end{lemma}

Hoeffding's inequality states that bounded centered random variables are sub-Gaussian with index proportional to the square of its oscillation:

\begin{proposition}
\label{hoeff}
Let $X$ be centered, with values in $[-b,a-b]$. We have that
\[
\| X \|_{\psi_2} \leq \tfrac{1}{2} a
\]
for every $\theta \in \bb R$.
\end{proposition}

Sub-Gaussian random variables behave like if they were Gaussian random variables, as long as moment bounds are concerned. In particular, we have \emph{Hanson-Wright inequality} (cf.~Lemma F.13 in \cite{JarMen}):

\begin{proposition}
\label{pa2}
Let $((X_i, Y_i); i \in I)$ be independent vectors in $\bb R^2$, and assume that $X_i$ and $Y_i$ are sub-Gaussian of indices $\sigma^2_i$ and $\tilde{\sigma}^2_i$. Let $g: I \times I \to \bb R$. We have that 
\[
\int \exp\Big\{ \gamma \sum_{i\neq j} g_{ij} X_i Y_j  \Big\} d \mu \leq 3
\]
whenever
\[
\gamma \leq \Big( 1024 \sum_{i \neq j} \sigma_i^2 \tilde{\sigma}^2_j g_{ij}^2 \Big)^{-1/2}.
\]
\end{proposition}

\thanks{{ \bf Acknowledgements:} J.A. and Y.X. thanks CAPES. M.J. has been funded by CNPq grant 201384/2020-5 and FAPERJ grant E-26/201.031/2022.}

\bibliographystyle{plain}


%
%
%
%
%
%
%
%
%
%
%

\end{document}